\documentstyle[12pt]{article}
\begin{document}
\newcommand{\C}{\it{I} \hspace{-2,8mm}C}

\def\L{{\cal L}}
\def\E{{\cal E}}
\def\K{{\cal K}}
\def\H{{\cal H}}
\def\U{{\cal U}}
\def\V{{\cal V}}
\def\R{{\cal R}}
\def\W{{\cal W}}
\def\Z{{\cal Z}}
\def\F{{\cal F}}
\def\M{{\cal M}}
\def\D{{\cal D}}
\def\G{{\cal G}}
\def\H{{\cal H}}
\def\A{{\cal A}}
\def\B{{\cal B}}
\def\diagram#1{\def\normalbaselines{\baselineskip=0pt\lineskip=10pt\lineskiplimit=1pt} \matrix{#1}}
\def\hfl#1#2{\smash{\mathop{\hbox to 8mm{\rightarrowfill}}\limits^{\scriptstyle#1}_{\scriptstyle#2}}}\def\vfl#1#2{\llap{$\scriptstyle #1$}\left\uparrow
\vbox to 6mm{}\right.\rlap{$\scriptstyle #2$}}
\def\antihfl#1#2{\smash{\mathop{\hbox to 8mm{\leftarrowfill}}\limits^{\scriptstyle#1}_{\scriptstyle#2}}}\def\antivfl#1#2{\llap{$\scriptstyle #1$}\left\downarrow
\vbox to 6mm{}\right.\rlap{$\scriptstyle #2$}}
\mathchardef\circ="220E
\newtheorem{defi}{D{\'e}finition}[section]
\newtheorem{prop}{Proposition}[section]
\newtheorem{theo}{Th{\'e}or{\`e}me}[section]
\newtheorem{lemme}{Lemme}[section]
\newtheorem{cor}{Corollaire}[section]
\def \lin {,\ldots ,}\def \diag#1{\def\normalbaselines{\baselineskip20pt\lineskip3pt \lineskiplimit3pt}\matrix{#1}}
\def \mapright#1{\smash{\mathop{\longrightarrow}\limits^{#1}}}
\def \mapdown#1{\Big\downarrow\rlap{$\vcenter{\hbox{$\scriptstyle#1$}}$}}
\def \mapup#1{\Big\uparrow\rlap{$\vcenter{\hbox{$\scriptstyle#1$}}$}}
\centerline{\bf Unitaire multiplicatif $K$-moyennable   }
\vskip4mm

\centerline{Mohamed Maghfoul}

\vskip6mm
{\bf R{\'e}sum{\'e}.}
\vskip2mm {\em Nous g{\'e}n{\'e}ralisons la th{\'e}orie de la $K$-moyennabilit{\'e} au cas d'un unitaire multiplicatif r{\'e}gulier $V$. Nous montrons que si $(H,V,U)$ est un syst\` eme de Kac  $K$-moyennable, alors pour toute  $S$-alg{\`e}bre $A$, les alg{\`e}bres $ A\times_{m}\hat S$ et $ A\times \hat S$ sont $KK$-{\'e}quivalentes o{\`u} $S$ est la $C^*$-alg{\`e}bre  de Hopf r{\'e}duite associ{\'e}e {\`a} $V$ .}
\vskip3mm
{\bf Introduction.}
\vskip3mm
Les $C^{*}$-alg{\`e}bres de Hopf peuvent {\^e}tre consid{\'e}r{\'e}es comme objets g{\'e}n{\'e}ralisant les groupes localement compacts. Elles ont {\'e}t{\'e} introduites dans le cadre de la recherche d'une cat{\'e}gorie contenant {\`a} la fois les groupes et leurs duaux et o{\`u} la dualit{\'e} de Pontrjagyn s'exprime ais{\'e}ment.
\vskip2mm
 Soient $H$ un espace de Hilbert et $V\in\L(H\otimes H)$ un op{\'e}rateur unitaire. On dit que $V$ est multiplicatif  s'il v{\'e}rifie, dans $\L(H\otimes H\otimes H)$,  la relation pentagonale 
$$V_{12}V_{13}V_{23} = V_{23}V_{12}.$$ 
Ces unitaires multiplicatifs introduits dans [2], fournissent un cadre g{\'e}n{\'e}ral contenant de nombreux exemples de groupes quantiques ( alg{\`e}bres de Kac, pseudo-groupes compacts matriciels de Woronovicz, groupes quantiques localement compacts,  groupes localement compacts etc ... )  o{\`u} la dualit{\'e} de Pontrjagyn a lieu.
\vskip2mm
Dans [1], S. Baaj et G. Skandalis ont associ{\'e} {\`a} toute paire $(A,B)$ de $C^{*}$-alg{\`e}bres munies d'une coaction d'une $C^{*}$-alg{\`e}bre de Hopf $S$ un groupe ab{\'e}lien not{\'e} $KK_{S}(A,B)$ de la th{\'e}orie de Kasparov {\'e}quivariante. Dans le cas o{\`u} $S = C_{0}(G)$ ($G$ est un groupe localement compact  et $C_{0}(G)$ est l'alg{\`e}bre des fonctions continues sur $G$ nulles {\`a} l'infini), $KK_{S}(A,B) = KK_{G}(A,B)$ [1] . \` A tout unitaire multiplicatif r{\'e}gulier sont associ{\'e}es quatre $C^*$-alg\` ebres de Hopf $ S, S_{p}, \hat{S}$ et $\hat{S_{p}}$. On dit que $V$ est $K$- moyennable s'il existe un $S$-module de Fredholm $(K,F,X)$ dans la classe de $1_{S}$ dans $ KK_{S}(\C,\C)$ et tel que la repr{\'e}sentation de $\hat{S_{p}}$ dans $K$ induite par la corepr{\'e}sentation $X$ de $S$ se factorise \` a travers la $C^*$-alg\` ebre r{\'e}duite dual $\hat{S}$. Dans le cas o{\`u} $V = V_{G}$: unitaire multiplicatif associ{\'e} au groupe $G$, la $K$-moyennabilit{\'e} de $V_{G}$ est {\'e}quivalente {\`a} celle de $G$ ([4], [5]). Un unitaire multiplicatif moyennable est $K$-moyennable. Nous montrons que si $(H,V,U)$ est un syst\` eme de Kac [2] avec $V$ est $K$-moyennable, alors pour toute $S$-alg{\`e}bre $A$, les alg{\`e}bres $ A\times_{m}\hat S$ et $ A\times \hat S$  [2] sont $KK$-{\'e}quivalentes. Si $V$ est discret (i.e. $\hat{S_{p}}$ est unif\` ere), la r{\'e}ciproque est vraie.
\vskip3mm
{\bf Remerciements.} Je tiens \`a remercier G. Skandalis pour les discussions enrichissantes que nous avons eues pendant l'\'elaboration de ce travail.
\vskip4mm 
\section{Rappels sur les unitaires multiplicatifs.}
 Soient $H$ un espace de Hilbert et $\L(H)$ l'alg{\`e}bre des op{\'e}rateurs born{\'e}s sur $H$.  Pour  $V\in\L(H\otimes H)$, nous d{\'e}signerons (dans $\L(H\otimes H\otimes H))$ par $V_{12}$, $V_{23}$ et $V_{13}$ les op{\'e}rateurs $V\otimes 1$, $1\otimes V$  et  $(\sigma\otimes 1)( 1\otimes V)(\sigma\otimes 1 )$ o{\`u}  $\sigma\in\L(H\otimes H)$ est la volte ( $\sigma(\xi\otimes\eta) = \eta\otimes\xi$ pour tous  $\xi,\; \eta\in H $).

\begin{defi}([2]) Soient $H$ un espace de Hilbert et $V\in\L(H\otimes H)$ un op{\'e}rateur unitaire. On dit que $V$ est multiplicatif s'il v{\'e}rifie la relation pentagonale $$V_{12}V_{13}V_{23} =V_{23}V_{12}.$$  
\end{defi}

 Par exemple, pour $G$ un groupe localement compact et $H = L^{2}(G)$, l'unitaire $V_{G}\in\L( L^{2}(G)\otimes L^{2}(G))$  tel que  $V_{G}(f)(s,t) = f(st,t)$ est multiplicatif. Pour d'autres exemples d'unitaires multiplicatifs voir[2].
\vskip3mm
{\`A} tout unitaire multiplicatif $V\in\L(H\otimes H)$, sont associ{\'e}s les deux sous-espaces de $\L(H)$ suivants:
$$A(V) = \{ L(\omega) = (\omega\otimes id)(V)/ \omega\in\L(H)_{\ast}\};$$
$$\hat{A}(V) = \{\varrho(\omega) = (id\otimes \omega)(V)/ \omega\in\L(H)_{\ast}\}.$$
On dit que $V$ est r{\'e}gulier si l'adh{\'e}rence normique de l'alg{\`e}bre $C(V) = \{(id\otimes \omega )(\sigma V)/ \omega\in\L(H)_{\ast} \}$ est $\K(H)$ ( alg{\`e}bre des op{\'e}rateurs compacts sur $H$). Notons que si $V$ est r{\'e}gulier, les adh{\'e}rences normique de $A(V)$ et $\hat{A}(V)$ qu'on note respectivement $S$ et $\hat{S}$ sont des $C^{*}$-alg{\`e}bres [2]. Munie du coproduit  $\delta(x) = V(x\otimes 1)V^{*}$, l'alg{\`e}bre $S$ de $V$ est une $C^{*}$-alg{\`e}bre de Hopf bisimplifiable. Munie du coproduit   $\hat{\delta} (x) = V^{*}(1\otimes x)V$, l'alg{\`e}bre $\hat{S}$ de $V$ est une $C^{*}$-alg{\`e}bre de Hopf bisimplifiable. L'alg{\`e}bre $S$ (resp. $\hat{S}$) est dite   alg{\`e}bre r{\'e}duite (resp. r{\'e}duite dual) de $V$.
\vskip2mm

\begin{defi} Soit $V\in\L(H\otimes H)$ un unitaire multiplicatif. Une repr{\'e}sentation de $V$ dans l'espace de Hilbert $K$ est un unitaire $X\in\L(K\otimes H)$ tel que $$X_{12}X_{13}V_{23} = V_{23}X_{12}.$$ 
Une corepr{\'e}sentation de $V$ dans un espace de Hilbert $K$ est la donn{\'e}e d'un unitaire $W\in\L(H\otimes K)$ v{\'e}rifiant $$V_{12}W_{13}W_{23} = W_{23}V_{12}.$$
 \end{defi}

 L'unitaire $V\in\L(H\otimes H)$ est une repr{\'e}sentation de $V$ dite repr{\'e}sentation r{\'e}guli{\`e}re. L'identit{\'e} de $H$ est appel{\'e}e  repr{\'e}sentation triviale de $V$. Si $(K,X)$ et $(K',X')$ sont deux repr{\'e}sentations de $V$, alors $(K\otimes K',X'_{23}X_{13})$ est une repr{\'e}sentation de $V$ dite produit tensoriel de $(K,X)$ par $(K',X')$. Pour toute repr{\'e}rentation (resp. corepr{\'e}sentation) $X$ de $V$ et pour tout $\omega\in\L(H)_{\ast}$ on pose $\varrho_{X}(\omega) = (id\otimes\omega)(X)$(resp. $L_{X}(\omega)= (\omega\otimes id)(X)$).

\begin{defi} 1) Soit $(S,\delta)$ une $C^{*}$-alg{\`e}bre de Hopf. Une corepr{\'e}sentation de $(S,\delta)$ dans l'espace de Hilbert ( ou $C^*$-module) $K$  est la donn{\'e}e d'un unitaire $U\in\L(K\otimes S)$ tel que $(id\otimes\delta)(U) = U_{12}U_{13}$ dans $\L(K\otimes S\otimes S)$.
 \vskip3mm 
2) Soit $A$ une $C^{*}$-alg{\`e}bre munie d'une coaction $\delta_{A}$ d'une $C^{*}$-alg{\`e}bre de Hopf $(S,\delta)$.  Une repr{\'e}sentation covariante de $(A,\delta_{A})$ est une paire $(\pi, U)$ o{\`u} $\pi$ est une repr{\'e}sentation   de $A$ et $U$ est une corepr{\'e}sentation de $S$ dans le m{\^e}me espace de Hilbert $K$ et telle que   
           $$(\pi\otimes id)\delta_{A}(a) = U(\pi(a)\otimes 1) U^{*}$$ 
pour tout $ a\in A$. 
\end{defi}

Pour $V$ un unitaire multiplicatif, on note $\hat{S}_{p}$ (resp. $S_{p}$ ) le s{\'e}par{\'e} compl{\'e}t{\'e} de l'espace $\L(H)_{\ast}$, muni du produit $(\omega\ast\omega')(T) = (\omega\otimes\omega')(V(T\otimes 1)V^{*})$ (resp.$(\omega\ast\omega')(T) = (\omega\otimes\omega')(V^*(1\otimes T)V))$ et de la semi-norme $\Vert \omega\Vert _{\hat{p}}= Sup \{\Vert \varrho_{X}(\omega)\Vert / \;   X \; \; repr\acute{e}sentation \; \; de \; \; V \}$ \hspace{3mm}        (resp. $\Vert   \omega\Vert _{p}= Sup\; \;  \{\Vert L_{X}(\omega)\Vert   /\; \; \-  X   \;\;   corepr\acute{e}sentation\;\; de\;\;  V \})$. Si $V$ est r{\'e}gulier,  $\hat{S}_{p}$ et $S_{p}$ sont des $C^*$-alg{\`e}bres de Hopf [2]. 
\vskip2mm
Il existe un unique unitaire $V'\in M(\hat{S}_{p}\otimes S)$ (resp. $V''\in M(\hat{S}\otimes S_{p}$)) tel que pour toute repr{\'e}sentation non d{\'e}g{\'e}n{\'e}r{\'e}e $\pi$ de $\hat{S}_{p}$ (resp. $S_{p}$), l'unitaire $X = (\pi\otimes L)(V')$ (resp. ($\varrho\otimes\pi)(V'')$) soit une repr{\'e}sentation (resp. corepr{\'e}sentation ) de $V$ avec  $\pi = \varrho_{X}$ (resp. $\pi = L_{X}$)( o{\`u} $M(A)$ d{\'e}signe l'alg{\`e}bre des multiplicateurs de la $C^*$-alg{\`e}bre $A$). La correspondance $X\rightarrow \varrho_{X}$ (resp. $X\rightarrow L_{X}$)  d{\'e}finit une bijection entre l'ensemble des repr{\'e}senatations (resp. Corepr{\'e}\-sentations) de $V$ et celui des repr{\'e}sentations non d{\'e}g{\'e}n{\'e}r{\'e}es de $\hat{S}_{p}$ (resp. $S_{p}$).
 \vskip2mm
Soient $V\in\L(H\otimes H)$ un unitaire multiplicatif r{\'e}gulier et $A$ une $C^{*}$-alg{\`e}bre munie d'une coaction $\delta_{A}$ de la $C^{*}$-alg{\`e}bre de Hop $(S,\delta)$ de $V$. Soit $(\pi,X)$ une repr{\'e}sentation covariante de $(A,\delta_{A})$ dans l'espace de Hilbert $K$. L'adh{\'e}rence normique dans $\L(K)$ de l'espace engendr{\'e} par les produits $\{\pi(a)\varrho_{X}(\omega), a\in A\; \; et\; \;\omega\in\L(H)_{\ast}\}$ est une $C^{*}$-  alg{\`e}bre. Le s{\'e}par{\'e} compl{\'e}t{\'e} de $A\otimes_{alg}\L(H)_{\ast}$ pour la semi-norme 
$$\Vert  \sum a_{i}\otimes \omega_{i}\Vert  = Sup\{\Vert \sum\pi(a_{i})\varrho_{X}(\omega_{i})\Vert \}$$ quand $(\pi,X)$ parcourt l'ensemble  des repr{\'e}sentations  covariantes  de $(A,\delta_{A})$  est une $C^{*}$-alg{\`e}bre not{\'e}e $A\times_{m}\hat S$ et appel{\'e}e produit  crois{\'e} maximal de $A$ par la coaction $\delta_{A}$ de $S$. Soit $\pi_{L} = (id\otimes L)\delta_{A}$ la repr{\'e}sentation de $A$ dans le $A$-module hilbertien $A\otimes H$. Alors $(\pi_{L}, V_{23})$ est une repr{\'e}sentation covariante de $(A,\delta_{A})$. L'adh{\'e}rence normique $A\times\hat S$ dans $\L(A\otimes H)$ de l'espace vectoriel engendr{\'e} par  $$\{\pi_{L}(a)(1\otimes\varrho(\omega)) / \; a\in A\;; \omega\in\L(H)_{\ast}\}$$ est une $C^{*}$-alg{\`e}bre dite produit crois{\'e} r{\'e}duit de $A$ par la coaction $\delta_{A}$   de $S$.

\vskip3mm

 Un unitaire multiplicatif $V\in\L(H\otimes H)$ est bir{\'e}gulier s'il est r{\'e}gulier et que l'adh{\'e}rencs normique dans $\L(H)$ de $\{(\omega  \otimes id)(\sigma V) /  \;  \omega\in \L(H)_{*} \}$ est $\K(H)$.

\begin{defi} Un syst\` eme de Kac est un triplet $(H,V,U)$ o\` u $H$ est un espace de Hilbert, $V\in\L(H\otimes H)$ est un unitaire multiplicatif bir{\'e}gulier et $U\in \L(H)$ est un unitaire tel que
\vskip2mm
$\bullet$  $U^2 = 1\; \; et \; \; (\sigma(1\otimes U)V)^3 = 1$
\vskip2mm
$\bullet$  Les unitaires $\hat{ V} = \sigma(U\otimes 1)V(U\otimes 1)\sigma\; \;  et\; \; \tilde{V} = (U\otimes U)\hat{ V}(U\otimes U)$ sont  multiplicatifs.
\end{defi} 

\section{$K$-moyennabilit{\'e}.}

Dans [1], S. Baaj et G. Skandalis ont associ{\'e} {\`a} toute paire $(A,B)$ de $C^{*}$-alg{\`e}bres munies d'une coaction d'une $C^{*}$-alg{\`e}bre de Hopf $S$ un groupe ab{\'e}lien not{\'e} $KK_{S}(A,B)$ de la th{\'e}orie de Kasparov {\'e}quivariante. Dans le cas o{\`u} $S = C_{0}(G)$ ($G$ est un groupe localement compact  et $C_{0}(G)$ est l'alg{\`e}bre des fonctions continues sur $G$ nulles {\`a} l'infini), $KK_{S}(A,B) = KK_{G}(A,B)$. Pour $ A = B = \C$, l'anneau (unif{\`e}re) $KK_{S}(\C,\C)$ est d{\'e}fini ainsi:
\vskip3mm

\begin{defi} Un $S$-module de Fredholm est un triplet $(K,F,X)$ o{\`u}:
\vskip2mm 
\noindent $\bullet$  $K$ est un espace de Hilbert s{\'e}parable $ Z/2 Z$- gradu{\'e}.
 \vskip2mm \noindent$\bullet$$F\in \L(K)$ un op{\'e}rateur de degr{\'e} 1 tel que  $1 - F^{*}F$  et  $1 - FF^{*}$ sont dans $\K(K).$
\vskip2mm
 \noindent $\bullet$  $X\in\L(K\otimes S)$ est une corepr{\'e}sentation de $S$ v{\'e}rifiant  $ (1\otimes x)[ F\otimes 1 - X(F\otimes 1)X^{*}]\in \K(K\otimes S)$  pour tout $x\in S.$
\end{defi}

Deux modules de fredholm $(K_{1},F_{1},X_{1})$ et $(K_{2},F_{2},X_{2})$ sont dits unitairement {\'e}quivalents s'il existe un unitaire $T\in \L(K_{1},K_{2})$ tel que $TF_{1} = F_{2}T$. On note $E_{S}(\C,\C)$l'ensemble des $S$-modules de Fredholm unitairement {\'e}quivalents.
\vskip2mm
On dit qu'un {\'e}l{\'e}ment $(K,F,X)\in E_{S}(\C,\C)$ est d{\'e}g{\'e}n{\'e}r{\'e} si $$1-FF^{*} = 1 - F^{*}F = (1\otimes x)(F\otimes 1 - X(F\otimes 1)X^{*}) = 0$$ pour tout $x\in S$.
\vskip3mm
\noindent Une homotopie est un {\'e}l{\'e}ment de  $E_{S}(\C,\C\otimes C[0,1])$. L'anneau ${\bf KK_{S}(\C,\C)}$ est l'ensemble des classes d'homotopie des {\'e}l{\'e}ments de $E_{S}(\C,\C)$. La multiplication est donn{\'e}e par le produit de Kasparov [1]. L'unit{\'e} {\'e}tant $(\C,0,id)$.  Pour une pr{\'e}sentation plus compl{\`e}te sur la th{\'e}orie de Kasparov {\'e}quivariante par rapport {\`a} une $C^*$-alg\` ebre de Hopf voir [1]. 
\vskip2mm
Soient $H$ un espace de Hilbert s{\'e}parable et $V\in \L(H\otimes H)$ un unitaire multiplicatif r{\'e}gulier. Soient $S$, $S_{p}$, $\hat{S}$ et $\hat{S_{p}}$ les $C^{*}$-alg{\`e}bres de Hopf associ{\'e}es {\`a} $V$. {\`A} toute corepr{\'e}sentation  $X \in \L(K\otimes S)$ de $S$ sur un espace de Hilbert $K$, correspond une repr{\'e}sentation (non d{\'e}g{\'e}n{\'e}r{\'e}e) de $\hat{S_{p}}$ dans $K$ et r{\'e}ciproquememt. Cette repr{\'e}sentation est donn{\'e}e par l'unitaire $Y = (id \otimes L)(X)\in\L(K\otimes H)$ [2]. On note $\rho:\hat{S_{p}}\longrightarrow \hat{S}$ l'application canonique.

\begin{defi} On dit que $V$ est $K$- moyennable s'il existe un $S$-module de Fredholm $(K,F,X)$ dans la classe de $1_{S}$ dans $ KK_{S}(\C,\C)$ et tel que la repr{\'e}sentation de $\hat{S_{p}}$ dans $K$ induite par la corepr{\'e}sentation  de $S$ se factorise \` a travers la $C^*$-alg\` ebre r{\'e}duite dual $\hat{S}$.
\end{defi}
{\bf Exemples.}
\vskip2mm 
1) Soient $G$ un groupe localement compact et $V_{G}$ l'unitaire multiplicatif asooci{\'e} ($V_{G}(f)(s,t) = f(st,t)$). Alors $V_{G}$ est $K$-moyennnable si, et seulemnt si, $G$ est $K$-moyennable.
\vskip2mm
2) Un unitaire multiplicatif r{\'e}gulier est dit moyennable, si l'application $\rho: \hat S_{p}\longrightarrow \hat S$ est isom{\'e}trique. Ceci {\'e}quivaut {\`a}: il existe une suite  de vecteurs  $\xi_{n}$ de $H$ de norme 1  v{\'e}rfiant $\Vert V(\xi_{n}\otimes \eta) - \xi_{n}\otimes\eta\Vert $ tend vers z{\'e}ro quand $n$ tend  vers l'infini et ceci pour tout $\eta\in H$; ou encore, la repr{\'e}sentation triviale $\tau$  est faiblement contenue dans la r{\'e}guli{\`e}re $\rho$ (cf.[2],[3]). Un unitaire multiplicatif moyennable est donc $K$-moyennable ( la $K$-moyennabilit{\'e} est donn{\'e}e par l'{\'e}l{\'e}ment trivial $(\C,0,Id)$ ). En particulier, un unitaires multiplicatif de type compact est  $K$-moyennable. \vskip3mm
{\bf Remarque.}
\vskip2mm
En rempla\c cant $S$ par $\hat{S}$ et $S_{p}$ par $\hat{S_{p}}$, on peut d{\'e}finir la $K$-comoyennabilit{\'e} d'un unitaire multiplicatif. Tous les r{\'e}sultats qui suivent ont lieu pour les unitaires multiplicatif $K$-comoyennables avec les remplacements cit{\'e}s pr{\'e}c{\'e}demment. 
  
\section{Applications.}

Soit $(H,V,U)$ un syst\` eme de Kac et soient $ S, S_{p},  \hat S$ et $ \hat S_{p}$ les  $C^{*}$-alg{\`e}bres de Hopf associ{\'e}es {\`a} $V$. Comme dans le cas d'un groupe[5], la d{\'e}monstration des lemmes suivants repose sur le fait que le produit tensoriel de toute
repr{\'e}sentation de l'unitaire multiplicatif $V$ par la repr{\'e}sentation r{\'e}guli\` ere est {\'e}quivalente \` a la repr{\'e}sentation r{\'e}guli{\`e}re [2]. Pour $B$ une  $C^{*}$-alg{\`e}bre quelconque, on note  $M(B)$ la $C^{*}$-alg{\`e}bre des multiplicateurs de $B$. Soient $A$ une $S$-alg{\`e}bre et $A\times_{m} \hat S$(resp. $A\times \hat S$) le produit crois{\'e} maximal (resp. produit crois{\'e} r{\'e}duit) de $A$ par $S$. Soit $V'\in M(\hat S_{p}\otimes S)$ l'unitaire associ{\'e} {\`a} la repr{\'e}sentation universelle de $\hat S_{p}$ [2]. Soit $(\pi,\U )$ la repr{\'e}sentation universelle covariante de $(A,\delta_{A})$ dans l'espace de Hilbert $\H$. Notons  $i_{A}: A\rightarrow M(A\times_{m}\hat S)$ et $i_{\hat S_{p}}:\hat S_{p}\longrightarrow M(A\times_{m}\hat S)$  les injections naturelles. On v{\'e}rifie facilement que la paire $(1\otimes\pi,\U_{23}V'_{13})$ d{\'e}finit une  repr{\'e}sentation cavariante de $(A,\delta_{A})$ dans $ M(\hat S_{p}\otimes A\times_{m}\hat S)$ et donc  un homomorphisme  $\Delta: A\times_{m} \hat S\longrightarrow M(\hat S_{p}\otimes A\times_{m}\hat S)$ ( $\otimes$ est le produit tensoriel minimal).
 
\begin{lemme} Il existe un homomorphisme   $\Delta_{1}:A\times \hat S\longrightarrow M(\hat S\otimes A\times_{m}\hat S)$ rendant commutatif le diagrammme suivant 
\vskip2mm
$$ \begin{array}{ccc} A\times_{m} \hat S&\mapright{\Delta}&M(\hat S_{p}\otimes A\times_{m}\hat S)\\
 \mapdown{\rho_{A}} & &\mapdown{\rho\otimes id}\\
 A\times \hat S&\mapright{\Delta_{1} }& M( \hat S\otimes A \times_{m}\hat S) \cr\\
\end{array}$$
\end{lemme}

 D{\'e}monstration. La repr{\'e}sentation $( \rho \otimes id)\circ\Delta$  de $A\times_{m} \hat S$  dans  $M(\hat S\otimes A\times_{m}\hat S)$ est donn{\'e}e par la repr{\'e}sentation covariante $(1\otimes i_{A},\U_{23}V_{13})$. Notons que $\hat S\otimes A \times_{m}\hat S$ se repr{\'e}sente fid\` element  dans  $\L( H\otimes \H)$. Mais  dans $ H\otimes \H$ les repr{\'e}sentations $\U_{23}V_{13}$  et  $V_{13}$ sont {\'e}quivalentes ( elles sont entrelac{\'e}es par l'unitaire $sXs^*$ o\` u $s:\H\otimes H\longrightarrow H\otimes \H$ est la volte). La repr{\'e}sentation  $( \rho \otimes id)\circ\Delta$  se factorise donc  \` a travers $A\times \hat S$. Ceci donne $\Delta_{1}$.

 \begin{lemme} Il existe un homomorphisme   $\Delta_{2}:A\times\hat S\longrightarrow M( \hat S_{p}\otimes A \times\hat S)$ rendant commutatif le diagramme suivant
 \vskip2mm
$$ \begin{array}{ccc} A\times_{m}\hat S&\mapright{\Delta}&M(\hat S_{p}\otimes A\times_{m}\hat S)\\
 \mapdown{\rho_{A}} & &\mapdown{ id\otimes\rho_{A}}\\
 A\times \hat S&\mapright{\Delta_{2} }& M( \hat S_{p}\otimes A \times\hat S) \cr\\
\end{array}$$
\end{lemme}
 D{\'e}monstration. Remarquons que la repr{\'e}sentation $(id\otimes\rho_{A})\circ\Delta$ est donn{\'e}e par la repr{\'e}srentation covariante $(1\otimes i_{A}, V_{23}V'_{13})$ o\` u $i_{A}$ est l'inclusion naturelle de $A$ dans $M( A \times\hat S)$. Soient $\varphi: A\longrightarrow \L(K)$ et $\psi: \hat S_{p} \longrightarrow \L(H_{p})$ des  repr{\'e}sentations fid\` eles. L'alg\` ebre $ A \times\hat S$ se repr{\'e}sente fid\` element dans  $\L(K\otimes H)$; et donc $ \hat S_{p}\otimes A\times\hat S$ se repr{\'e}sente fid\` element dans  $\L(H_{p}\otimes K\otimes H)$. Soit $\tilde{\psi}$ une telle repr{\'e}sentation. Alors $(\tilde{\psi}\otimes id )(V_{23}V'_{13)}) = V_{34}W_{14}$ dans $\L(H_{p}\otimes K\otimes H\otimes H)$ o\` u $W = (\psi\otimes id )(V')$. Comme $W$ est une repr{\'e}sentation de $V$ et qu'on a \` a faire \` a un syst\` eme de Kac, la  repr{\'e}sentation  $V_{13}W_{23} $ est {\'e}quivalente \` a $V_{13}$  dans $\L(H\otimes H_{p}\otimes H)$  ( voir d{\'e}monstration de la proposition A.10 de [2]), par suite  $V_{23}W_{13} $ est {\'e}quivalente \` a $V_{23}$ dans $\L(H_{p}\otimes H\otimes H)$ . Donc  $V_{34}W_{14} $ est {\'e}quivalente \` a $V_{34}$. Par cons{\'e}quent $(id\otimes\rho_{A})\circ\Delta$ se factorise \` a travers $ A\times \hat S$. Ceci donne $\Delta_{2}$.       \vskip3mm
 Les diagrammes des deux lemmes commutent par construction.

\begin{theo} Soient $ (H,V,U)$ un syst\` eme de Kac, $ S$ la $C^*$-alg\` ebre r{\'e}duite associ{\'e}e \` a $V$ et  $A$ une $S$-alg{\`e}bre. Si $V$ est $K$-moyennable, alors l'application $\rho_{A}: A\times_{m}\hat S \longrightarrow A\times \hat S$ est inversible en $KK$-th{\'e}orie.\end{theo}

D{\'e}monstration.  Soient $i_{A}: A\rightarrow M(A\times_{m}\hat S)$ et $i_{\hat S_{p}}:\hat S_{p}\longrightarrow M(A\times_{m}\hat S)$  les injections naturelles. Soit $V'\in M(\hat S_{p}\otimes S)$ l'unitaire associ{\'e} {\`a} la  repr{\'e}sentation universelle de $\hat S_{p}$ [2]. On pose $\W =(i _{\hat S_{p}}\otimes id_{S})(V') \in M ( A\times_{m}\hat S\otimes S)$ . Soit   $\Phi_{A}: KK_{S}(\C,\C)\longrightarrow  KK(A\times_{m}\hat S,A\times_{m}\hat S)$   le morphisme de
 Kasparov[2]. L'image par  $\Phi_{A}$ de l'{\'e}l{\'e}ment $(K,F,X)$ r{\'e}alisant  la $K$-moyennabilit{\'e} de  $V$ est de la forme $(K\otimes A\times_{m}\hat S, F\otimes 1)$. L'action de $ A\times_{m}\hat S$ dans $(K\otimes A\times_{m}\hat S)$ est donn\' ee par la repr\' esentation covariante   $ (i_{A},\W_{23}X_{13})$( l'alg\` ebre $A$ agit comme multiplicateur  sur le deuxi\` eme facteur de $K\otimes A\times_{m}\hat S)$.  Cette repr{\'e}sentation provient d'une  repr{\'e}sentation de $\hat S_{p}\otimes A\times_{m} \hat S$ qui, par hypoth{\`e}se, se factorise {\`a}  travers $\hat S\otimes A\times_{m}\hat S$ et par cons{\'e}quent (lemme 3.1) {\`a} travers $A\times\hat S$. On obtient donc ( puisque  $\Phi_{A}(1) =  1$) un {\'e}l{\'e}ment $\omega_{A}\in KK(A\times\hat S,A\times_{m}\hat S)$ tel que $\lambda_{A}\otimes_{A\times\hat S}\omega_{A} = 1$.  
\vskip2mm
\noindent Montrons que $\omega_{A}\otimes_{A\times_{m}\hat S}\lambda_{A} = 1$. Notons que l'{\'e}l{\'e}ment $\omega_{A}\otimes_{A\times_{m}\hat S}\lambda_{A}$ est de la forme $ (K\otimes A\times\hat S, F\otimes 1)$. L'action de $\hat S_{p}\otimes A\times_{m}\hat S$ dans  $K\otimes A\times\hat S$ se factorise {\`a} travers $\hat S_{p}\otimes A\times\hat S$. Soit  $x_{t} = (K_{t},F_{t},X_{t})$ une homotopie dans $KK_{S}(\C,\C)$ entre $x_{0} = (K,F,X)$ et $x_{1}= 1_{S}$. Pour chaque $t$, nous avons une action de  $\hat S_{p}\otimes A\times\hat S$ sur  $K_{t}\otimes A\times\hat S$. Ceci  et le lemme (3.2) montrent que l'action de $A\times_{m}\hat S$  se factorise {\`a} travers $A\times\hat S$. Nous en d{\'e}duisons que $(K_{t}\otimes A\times \hat S, F_{t}\otimes 1)$ r{\'e}alise une homotopie entre $\omega_{A}\otimes_{A\times_{m}\hat S}\lambda_{A}$ et $1_{A\times\hat S}$.
\vskip3mm
 Rappelons (voir [7]) qu'une  $S$-alg{\`e}bre $A$ est dite nucl{\'e}aire en $KK_{S}$-th{\'e}orie si $1_{A}\in KK_{S}(A,A)$ est repr{\'e}sent{\'e} par un $A,A$-bimodule nucl{\'e}aire.

\begin{prop}Soit $(H,V,U)$ un syst\` eme de Kac. Si $V$ est $K$-moyennable et $A$ est nucl{\'e}aire en $KK_{S}$-th{\'e}orie, alors les alg{\`e}bres $A\times\hat S$ et $A\times_{m}\hat S$ sont $K$-nucl{\'e}aires.
 \end{prop}

 D{\'e}monstration. Il suffit de montrer que $A\times_{m}\hat S$ est $K$-nucl{\'e}aire car  $A\times\hat S$ et  $A\times_{m}\hat S$ sont $KK$-{\'e}quivalentes (th{\'e}or{\`e}me 3.1). Soient $(H_{1},F_{1})$ l'{\'e}l{\'e}ment r{\'e}alisant la $K$-moyennabilit{\'e} de $V$ et $(E_{2},F_{2})$ l'{\'e}l{\'e}ment r{\'e}alisant la nucl{\'e}arit{\'e}  en $KK_{S}$-th{\'e}orie de  $A$. Soit $D$ une $C^{*}$-alg{\`e}bre sur laquelle $S$ agit trivialement. Il existe $F$ (une $ F_{2}$-connection ) tel que $$ ((H_{1}\otimes E_{2})\times_{m}\hat S, F) = 1_{A\times_{m}\hat S}.$$ D'autre part  $$((H_{1}\otimes E_{2})\times_{m}\hat S)\otimes_{max}D = (H_{1}\otimes A\otimes_{min}D)\times_{m}\hat S \otimes_{(A\otimes_{min}D)\times_{m}\hat S}(E_{2}\otimes_{max}D)\times_{m}\hat S.$$ D'apr{\`e}s le th{\'e}or{\`e}me (3.1), $(H_{1}\otimes A\otimes_{min}D)\times_{m}\hat S$ est un $((A\otimes_{min}D)\times\hat S,(A\otimes_{min}D)\times_{m}\hat S)$ bimodule. Comme $(A\otimes_{min}D)\times\hat S$ est un quotient de $A\times_{m}\hat S\otimes_{min}D$, alors  $(H_{1}\otimes E_{2})\times_{m}\hat S$ est nucl{\'e}aire.
\vskip3mm 
\begin{cor} Soit $(H,V,U)$  un syst\` eme de Kac avec  $V$  discret. Les propri{\'e}t{\'e}s suivantes sont {\'e}quivalentes.
\vskip2mm
 i) $V$ est $K$-moyennble.
\vskip2mm
ii) Pour toute $S$-alg{\`e}bre $A$, l'application $\rho_{A}: A\times_{m}\hat S\longrightarrow A\times \hat S$ est inversible en $KK$-th{\'e}orie.
\vskip2mm
iii) L'application $\rho^{*}: KK(\hat S, \C)\longrightarrow  KK(\hat S_{p}, \C)$ est un isomorphisme.
\end{cor}

D{\'e}monstration. $i)\Longrightarrow ii)$ (th{\'e}or{\`e}me 3.1) et $ii)\Longrightarrow iii)$ {\'e}vident.
\noindent $iii)\Longrightarrow i)$. Parce que $V$ est discret, on a $KK_{S}(\C,\C) = KK(\hat S_{p},\C)$. Il existe donc $\;  x\in KK(\hat S, \C)$ tel que $\rho^{*}(x) = 1 $ dans $KK_{S}(\C,\C)$. La repr{\'e}sentation de $\hat S_{p}$ dans l'espace de Hilbert du $S$-module de Fredholm $\rho^{*}(x)$ provient  d'une repr{\'e}sentation de $\hat S$, elle est donc faiblement contenue dans la repr{\'e}sentation  r{\'e}guli{\`e}re $\rho$.

\vskip6mm
{\em Universit{\'e} Ibn Tofail

D{\'e}partement de Math{\'e}matiques

BP 133

K{\'e}nitra, Maroc}

e-mail: mmaghfoul@mailcity.com

\end{document}